\documentclass[11pt,leqno]{article}

\usepackage{amsmath,amsthm}
\usepackage {latexsym}
\usepackage{amssymb}

\newcommand{\bse}{\begin{equation}}

\newcommand{\bear}{\begin{eqnarray}}
\newcommand{\eear}{\end{eqnarray}}

\newcommand{\supp}{\mbox{\rm supp}}

\newcommand{\half}{\frac{1}{2}}

\newcommand{\eps}{{\varepsilon}}
\newcommand{\R}{{\mathbb R}}

\newcommand{\Compl}{{\mathbb C}}

\newcommand{\les}{\lesssim}

\newcommand{\Laplace}{\triangle}

\newcommand{\kato}{{\mathcal K}}

\def\Lap{\Delta}
\def\nn{\nonumber}

\newtheorem{theorem}{Theorem}
\newtheorem{lemma}[theorem]{Lemma}
\newtheorem{defi}[theorem]{Definition}
\newtheorem{cor}[theorem]{Corollary}
\newtheorem{prop}[theorem]{Proposition}
\newtheorem{proposition}[theorem]{Proposition}

\theoremstyle{remark}

\def\Pac{P_{a.c.}}
\def\la{\langle}
\def\ra{\rangle}

\def\norm[#1][#2]{\|#1\|_{#2}}

\def\y{{\bf y}}
\def\x{{\bf x}}

\def\disp{\displaystyle}
\def\norm[#1][#2]{\Vert #1 \Vert_{#2}}

\begin{document}

\title{Dispersive estimates for the three-dimensional Schr\"odinger equation
 with rough potentials}
\date{}

\author{M.\ Goldberg}

\maketitle

\begin{abstract}
      The three-dimensional Schr\"odinger propogator $e^{itH}$,
$H = -\triangle + V$, is a bounded map from $L^1$ to $L^\infty$
with norm controlled by $|t|^{-3/2}$ provided the potential
satisfies two conditions:  An integrability condition limiting the
singularities and decay of $V$, and a zero-energy spectral condition
on $H$.  This is shown by expressing the spectral measure of $H$ in
terms of its resolvents and proving a family of $L^p$ mapping estimates
for the resolvents.  Previous results in this direction had required
$V$ to satisfy explicit pointwise bounds.
\end{abstract}

\section{Introduction}

In this paper we consider dispersive estimates for the the time evolution
operator $e^{itH}P_{ac}(H)$, where $H = -\Lap + V$ in $\R^3$ and $P_{ac}(H)$ 
is the projection onto the absolutely continuous subspace of $H$.  Our goal is
to assume as little as possible on the potential $V=V(x)$ in terms of decay
or regularity. More precisely, we prove the following theorem.

\begin{theorem}
\label{thm:dispersive}
Let $V \in L^{\frac32(1+\eps)}(\R^3) \cap L^1(\R^3)$.
Assume also that zero is neither an eigenvalue nor a resonance of $H=-\Laplace+V$. 
Then
\begin{equation}
\label{eq:dis_3} \big\|e^{itH} P_{ac}(H)\big\|_{1\to\infty} \les 
|t|^{-\frac32}. 
\end{equation}
\end{theorem}

See Section~3 for a discussion of resonances.  With this assumption the 
spectrum is known to be purely absolutely continuous on $[0,\infty)$, 
see \cite{GS2} for details.
 
Previous results in this direction have generally required pointwise decay of
$V$.  Journ\'e, Soffer and Sogge \cite{JSS} proved a version of 
Theorem~\ref{thm:dispersive} under the pointwise bound 
$|V(x)| \le C(1+|x|)^{-\beta},\, \beta>7$, and also some regularity 
assumptions including $\hat{V} \in L^1$.
Yajima~\cite{Y1} reduced the decay
hypothesis to $\beta > 5$ and proved that the wave operators are bounded
on $L^p(\R^3)$ for all $1 \le p \le \infty$.  The dispersive estimate 
follows from this result.  Finally, Goldberg and Schlag~\cite{GS1}
established the dispersive estimate provided $\beta >3$.  In all these works
the assumption is made that zero energy is neither an eigenvalue nor a
resonance.


The exposition in this paper roughly follows~\cite{GS1}, with two significant
refinements.  First, the distinction which was previously drawn between high 
and low energies is now removed.  Second, the limiting absorption principle
of Agmon~\cite{agmon}, which concerns the action of resolvents on weighted 
$L^2$, is replaced with unweighted $L^p$ estimates as in \cite{GS2}.
The ability to work with potentials that satisfy $L^p$ conditions (but not 
necessarily any pointwise bounds) depends in turn on a unique continuation 
result due to Ionescu and Jerison~\cite{IonJer}.  For reference we present
the statement here.
\begin{theorem} \label{thm:IJ}
Let $V\in L^{\frac32}(\R^3)$. Suppose $u\in W^{1,2}_{\rm loc}(\R^3)$ satisfies
$(-\Laplace +V)u=\lambda^2 u$ where $\lambda\ne0$ in the sense of
distributions. If, moreover,
$\|(1+|x|)^{\delta-\half}u\|_2<\infty$ for some $\delta>0$, then $u\equiv 0$.
\end{theorem}

In terms of local regularity, Theorem~\ref{thm:dispersive} appears to be
nearly optimal.  There exist compactly supported potentials 
$V \in L^{3/2}_{\rm weak}$ for which $-\Lap + V$ admits bound states with 
positive energy~\cite{KocTat}.  On the other hand, while the assumption
$V \in L^1(\R^3)$ corresponds to (radial) pointwise decay on the order of 
$|V(x)| \le C(1+|x|)^{-3-\eps}$, it is reasonable to expect dispersive 
behavior to persist even with weaker decay hypotheses on $V$.  This is already
shown in~\cite{RS} for small potentials in the Kato class, which includes
all $V\in L^{\frac32 + \eps} \cap L^{\frac32 - \eps}$ of small norm.

\subsection{Resolvent Identities}

Let $H=-\Lap+V$ in $\R^3$ and define the resolvents
$R_0(z) :=(-\Lap-z)^{-1}$ and $R_V(z):= (H-z)^{-1}$.
For $z \in \Compl \setminus \R^+$, the operator $R_0(z)$ can be realized as an
integral operator with the kernel
\[R_0(z)(x,y) = \frac{e^{i\sqrt{z}|x-y|}}{4\pi|x-y|} \]
where $\sqrt{z}$ is taken to have positive imaginary part.
While $R_V(z)$ does not possess an explicit representation of this form,
it can be expressed in terms of $R_0(z)$ via the identities
\begin{equation} \begin{aligned} \label{eq:resident}
R_V(z) &= (I + R_0(z)V)^{-1}R_0(z) = R_0(z)(I + VR_0(z))^{-1} \\
R_V(z) &= R_0(z) -R_0(z)VR_V(z)    = R_0(z) - R_V(z)VR_0(z)
\end{aligned}
\end{equation}
In the case where $z = \lambda \in \R^+$, one is led to consider limits
of the form $R_0(\lambda \pm i0) := \lim_{\eps\downarrow 0} 
R_0(\lambda\pm i\eps)$.  The choice of sign determines which branch of the
square-root function is selected in the formula above, therefore the two
continuations do not agree with one another.  
 For convenience we will adopt a shorthand
notation for dealing with resolvents along the positive real axis, namely
\begin{align*}
R_0^\pm(\lambda) &:= R_0(\lambda\pm i0) \\
R_V^\pm(\lambda) &:= R_V(\lambda\pm i0)
\end{align*}
Note that $R_0^-(\lambda)$ is the formal adjoint of $R_0^+(\lambda)$,
and a similar relationship holds for $R_V^\pm(\lambda^2)$.
The discrepancy between $R_0^+(\lambda)$ and $R_0^-(\lambda)$
characterizes the absolutely continuous part of the spectral measure
of $H$, denoted here by $E_{ac}(d\lambda)$, by means of the Stone formula
\begin{equation}
\label{eq:stone}
\langle E_{ac}(d\lambda) f,g \rangle
= \frac{1}{2\pi i}\big\langle [R_V^+(\lambda)-R_V^-(\lambda)]f,g \big\rangle 
\,d\lambda.
\end{equation}

Let $\chi$ be a smooth, even, cut-off function
on the line that is equal to one on a neighborhood of the origin. 
In order to prove Theorem~\ref{thm:dispersive}
it will suffice to show that
\begin{align} 
\sup_{L\ge1}\Big|\Big \la e^{itH} \chi(\sqrt{H}/L)\Pac f,g \Big\ra\Big| &=  \sup_{L\ge1}\Big|\int_0^\infty 
e^{it\lambda^2}\lambda\, \chi(\lambda/L) \Big \la [R_V^+(\lambda^2)-
R_V^-(\lambda^2)]f,g \Big\ra
\, \frac{d\lambda }{\pi i}\Big| \label{eq:spec_theo} \\
&\les |t|^{-\frac32}\|f\|_1\|g\|_1. \nn
\end{align}
The first equality is precisely \eqref{eq:stone}, and we have also made the
change of variable $\lambda \mapsto \lambda^2$.

Our approach roughly parallels the one found in \cite{GS1}, with two main 
differences.  The first is that norms will be estimated in a variety of $L^p$
spaces in addition to the more typical weighted $L^2$.  The second
is that low and high energies will not require a separate calculation.  
There is still a distinction to be noted between the two cases, however.
The limiting absorption principle is used to establish decay as 
$\lambda\to\infty$, whereas boundedness at low energies follows from a 
Fredholm alternative argument.  This requires assuming that zero is neither an
eigenvalue nor a resonance.
 
\subsection{Initial terms of the Born series}
Iterating the resolvent identity \eqref{eq:resident} a total of $m+2$
times yields the finite Born series
\begin{align}
R_V^\pm(\lambda^2) &=  \sum_{k=0}^{m+1} R_0^\pm(\lambda^2)
(-VR_0^\pm(\lambda^2))^k \nn \\
& +R_0^\pm(\lambda^2)V R_V^\pm(\lambda^2)(VR_0^\pm(\lambda^2))^{m+1}.
\label{eq:res_ident}
\end{align}
Here $m$ is any positive integer.  This expansion is then inserted into
the integral in \eqref{eq:spec_theo}. 
The first $m+2$ terms which do not contain the resolvent $R_V$ 
are treated as in~\cite{RS}, Section~2, which only requires that
\begin{equation}
\label{eq:katonorm}
 \|V\|_{\kato}:=\sup_{x\in\R^3} \int \frac{|V(y)|}{|x-y|}\, dy < \infty. 
\end{equation}
In particular, if $V \in L^{\frac32+\eps} \cap L^{\frac32-\eps}$, 
then this condition is satisfied by dividing $\R^3$ into the regions $|x-y|<1$
and $|x-y| \ge 1$.

For the convenience of the reader we recall the relevant arguments 
from~\cite{RS}.  When the Born series \eqref{eq:res_ident} is substituted into
\eqref{eq:spec_theo}, the contricuation from the $k^{\rm th}$ term is equal to
\begin{equation*}   \int_0^\infty e^{it\lambda^2}\lambda\;\psi(\lambda/L) \,
  \big\la \big[R_0^+(\lambda^2)(VR_0^+(\lambda^2))^k - R_0^-(\lambda^2)
    (VR_0^-(\lambda^2))^k\big]\,f,g \big\rangle \, d\lambda
\end{equation*}
which is controlled by 
\begin{align}
& \sup_{L\ge1} \Bigl| \int_0^\infty e^{it\lambda^2}\lambda\;\psi(\lambda/L) \,
 \Im \langle R_0^+(\lambda^2)(VR_0^+(\lambda^2))^k\,f,g \rangle \,
d\lambda \Bigr| \nonumber\\
& \le\int_{\R^6} |f(x_0)||g(x_{k+1})|
\int_{\R^{3k}} \frac{\prod_{j=1}^k |V(x_j)|}{\prod_{j=0}^k 4\pi |x_j-x_{j+1}|}\cdot\nonumber\\
& \qquad\qquad\qquad \cdot \sup_{L\ge1} \Bigl| \int_0^\infty e^{it\lambda^2}
\lambda\; \psi(\lambda/L)\, \sin\Bigl(\lambda\sum_{\ell=0}^k 
|x_\ell-x_{\ell+1}|\Bigr)\, d\lambda \Bigr|
\; d(x_1,\ldots,x_k)\,dx_0\,dx_{k+1} \label{eq:gross}\\
&\le Ct^{-\frac32}  \int_{\R^6} |f(x_0)||g(x_{k+1})| 
\int_{\R^{3k}} \frac{\prod_{j=1}^k |V(x_j)|}{(4\pi)^{k+1}\prod_{j=0}^k|x_j-x_{j+1}|}\sum_{\ell=0}^k |x_\ell-x_{\ell+1}|
\; d(x_1,\ldots,x_k)\;dx_0\,dx_{k+1} \label{eq:lem1} \\
\end{align}
which in turn is controlled by
\begin{align}
&\le Ct^{-\frac32}  \int_{\R^6} |f(x_0)||g(x_{k+1})|\; 
(k+1) (\|V\|_{\kato}/4\pi)^k \;dx_0\,dx_{k+1} \label{eq:lem2} \\
&\le C_k\,t^{-\frac32} \|f\|_1\|g\|_1. \nonumber
\end{align}
In order to pass to \eqref{eq:gross} one uses the explicit representation 
of the kernel of 
$R_0^+(\lambda^2)(x,y)=\frac{e^{ i\lambda|\x-\y|}}{4\pi|\x-\y|}$, 
which leads to a $k$-fold integral. The inequalities~\eqref{eq:lem1} 
and~\eqref{eq:lem2} are obtained by means of the following two lemmas 
from~\cite{RS}, which we reproduce here without proof.
They may be regarded as exercises in the use of stationary phase and
Fubini's Theorem, respectively. 

\begin{lemma}
\label{lem:statphas}
Let $\psi$ be a smooth, even  bump function
with $\psi(\lambda)=1$ for $-1\le\lambda\le 1$ and $\supp(\psi)\subset[-2,2]$. 
Then for all $t\ge1$ and any real~$a$,
\begin{equation}
\label{eq:decay}
\sup_{L\ge 1}\Bigl| \int_0^\infty e^{it\lambda} \sin(a\sqrt{\lambda})\,
\psi\Bigl(\frac{\sqrt\lambda}{L}\Bigr)\, d\lambda\Bigr| \le
C \,t^{-\frac32}\,|a|
\end{equation}
where $C$ only depends on $\psi$ and $\chi$.
\end{lemma}

\begin{lemma}
\label{lem:iter} 
For any positive integer $k$ and $V$ as in \eqref{eq:katonorm}
\begin{equation}
\nn
\sup_{x_0,x_{k+1}\in\R^3}\int_{\R^{3k}} \frac{\prod_{j=1}^k |V(x_j)|}{\prod_{j=0}^k|x_j-x_{j+1}|}\sum_{\ell=0}^k |x_\ell-x_{\ell+1}|\; dx_1\ldots\,dx_k \le (k+1) \|V\|_{\kato}^k.
\end{equation}
\end{lemma}

\section{Estimates on the free resolvent}

We now turn to the term in the Born series~\eqref{eq:res_ident} containing the
perturbed resolvent~$R_V$.  The following propositions establish a 
family of mapping estimates for the free resolvent.
\begin{proposition} \label{prop:limabs}
For each exponent $1 < p \le \frac43$, there exist constants $C_p < \infty$
such that
$$\norm[R_0^\pm(\lambda^2)f][L^{3p}] \le C_p \lambda^{-2+2/p}
\norm[f][L^{p}]$$
For each exponent $\frac43 \le p < \frac32$, there exist constants 
$C_p < \infty$ such that
$$\norm[R_0^\mp(\lambda^2)f][L^{p*}] \le C_p\lambda^{4-6/p}
\norm[f][L^p] \qquad {\rm where}\ \frac1{p*} = \frac3p - 2$$
\end{proposition}
\begin{proof}
The case $p = \frac43$ is proven as a special case of theorem~2.3 in 
\cite{KRS}.  It is clear from fractional integration that
$R_0^\pm(\lambda^2)$ maps $L^1(\R^3)$ to weak-$L^3(\R^3)$ uniformly in 
$\lambda$, using the definition
$$\norm[f][L^3_{\rm weak}(\R^3)] = \sup_{A\subset \R^3, |A| < \infty}
  |A|^{-2/3}\int_A |f(x)|\,dx$$
which is equivalent to the usual weak-$L^3$ ``norm'' and also satisfies a
triangle inequality, see Lieb, Loss~\cite{LieLos}, Chapter 4.3.
  The cases $1 < p < \frac43$ and $\frac43 < p < \frac32$ follow
by Marcinkiewicz interpolation and duality, respectively.
\end{proof}
\begin{proposition} \label{prop:freeRV}
Suppose $V \in L^{\frac32(1+\eps)}(\R^3) \cap L^{\frac32(1-\eps)}(\R^3)$.
Then there exists a constant $C_\eps < \infty$ such that
\begin{equation} \label{eq:freeVR}
\norm[VR_0^\pm(\lambda^2)f][L^p]\le C_\epsilon(1+|\lambda|)^{-2+2/p} \norm[V][]
\norm[f][L^p] \quad {\it for\ all\ exponents\ }\ 1\le p \le 1+\eps.
\end{equation}
The dual operators satisfy the related bound
\begin{equation} \tag{\ref{eq:freeVR}'} \label{eq:freeRV}
\norm[R_0^\pm(\lambda^2)Vf][L^p] \le C_\eps (1+|\lambda|)^{-2/p} \norm[V][]
\norm[f][L^p] \quad {\it for\ all\ exponents\ }\ 
\frac{1+\eps}{\eps} \le p \le \infty.
\end{equation}
\end{proposition}
In the above statement $\norm[V][]$ is understood to be the larger of
$\norm[V][L^{\frac32(1+\eps)}]$ and $\norm[V][L^{\frac32(1-\eps)}]$.

\begin{proof}
In the case $p=1$, $(VR_0^\pm(\lambda^2))$ has an operator bound of precisely
$(4\pi)^{-1}\norm[V][\kato]$, which is controlled by $\norm[V][]$.  
The case $p=1+\epsilon$, $|\lambda| > 1$, is a corollary of the preceding
proposition, using the fact that $V \in L^{\frac{3p}2}$.  For $|\lambda|\le 1$,
a uniform bound is obtained by comparing $R_0(\lambda^2)$ to fractional
integration and observing that $V \in L^{\frac32}$.
The intermediate cases $1<p<1+\eps$ follow by interpolation.
\end{proof}

It is slightly inconvenient that $\norm[VR_0^\pm(\lambda^2)][1\to 1]$ does not 
decay in the limit $|\lambda|\to \infty$.  If this map is iterated several
times, however, we may use the fact that $(VR_0^\pm(\lambda^2))$ maps 
$L^1(\R^3)$ to $L^{2/(2-\eps)}(\R^3)$ and vice versa to apply the bound in 
\eqref{eq:freeVR} with $p = \frac2{2-\eps}$.  The resulting mapping estimates
will be needed in section~\ref{sec:calculations}.
\begin{cor}
Let $V \in L^{\frac32(1+\eps)}(\R^3)\cap L^{\frac32(1-\eps)}(\R^3)$. Then
\begin{equation} \label{eq:(RV)^k}
\begin{aligned}
\norm[(VR_0^\pm(\lambda^2))^{k+2}f][L^1] &\le C^k(1+|\lambda|)^{k\eps}
\norm[f][L^1] \qquad {\it and} \\
\norm[(R_0^\pm(\lambda^2)V)^{k+2}f][L^\infty] &\le C^k(1+|\lambda|)^{k\eps}
\norm[f][L^\infty]
\end{aligned} \end{equation}
\end{cor}

We now consider the action of the free resolvent on weighted $L^p$ spaces.
Let $L^{p,\sigma}(\R^3)$ be the Banach space determined by the norm
$$\norm[f][L^{p,\sigma}] = \norm[(1+|\cdot|)^{\sigma}f][L^p], \quad 
  1\le p\le\infty, \sigma \in \R$$

\begin{lemma} \label{lem:weightedRV}
Suppose $V\in L^{\frac32(1+\eps)}(\R^3) \cap L^1(\R^3)$, and let $p$ be
any exponent in the range $\frac{1+\eps}{\eps} \le p \le \infty$. 
The operator $R_0^\pm(\lambda^2)V$ is a bounded linear map on
$L^{p,-1}$, and its norm is controlled by $\norm[V][]$.

Furthermore, $R_0^\pm(\lambda^2)V$ may be written as 
a sum of linear maps $T_1$ and $T_2$ satisfying the following estimates:
\begin{align} \label{eq:weightedRV}
&\norm[T_1f][L^{p,-1}] \les (1+|\lambda|)^{-1/p}\norm[V][]\norm[f][L^{p,-1}] \\
&\norm[T_2f][L^\infty]\  \les \norm[V][] \norm[f][L^{p,-1}
  (\{|x|>\lambda^{1/p}\})]
\end{align}
The constant of similarity depends on $\eps > 0$ but not on
the specific choice of $p$.
\end{lemma}
\begin{proof}
For each $k = 0,1,2,\ldots$, let 
$D_k = \{x\in\R^3: |x| < \lambda^{1/p}2^{k+1}\}$.
We define $T_1$ and $T_2$ in the following manner: In the annulus 
$A_k = \{x: 2^{k-1} \le |x| < 2^k\}$ (or the unit ball $A_0 = \{|x| < 1\}$), 
let
$$ \begin{aligned}
T_1f(x) &= R_0^\pm(\lambda^2)V\chi_{D_k}f(x) \\
T_2f(x) &= R_0^\pm(\lambda^2)V\chi_{D_k^c}f(x)
\end{aligned} $$

The estimate for $T_2f$ is immediate.  Since $V \in L^{p'}$, by H\"older's
inequality $Vf \in L^{1,-1}$.  The construction of $D_k$ ensures that 
$|y-x| > (1 +|y|)/3$ for any $x\in A_k$, $y \in D_k^c$.  Thus
$$|T_2f(x)|\ < \ \frac3{4\pi} \int_{D_k^c} \frac{|V(y)f(y)|}{1+|y|}\,dy
\ < \ \frac3{4\pi} \norm[Vf][L^{1,-1}(D_0^c)]
\ \les \ \norm[V][]\norm[f][L^{p,-1}(D_0^c)]$$
It should be noted that $L^\infty(\R^3)$ has a natural embedding into
$L^{p,-1}(\R^3)$ for any $p > 3$.

To control $T_1f$, we first consider its restriction to each annulus $A_k$.
Proposition \ref{prop:freeRV} states that
$\norm[T_1f][L^p(A_k)] \les \norm[R_0^\pm(\lambda^2)V\chi_{D_k}f][L^p] \les
(1+|\lambda|)^{-2/p}\norm[V][]\norm[f][L^p(D_k)]$.  
The norm of $T_1f$, as measured in the space $L^{p,-1}(\R^3)$,is recovered by 
summing over $k$.
$$ \begin{aligned}
\norm[T_1f][L^{p,-1}]^p \ &\sim \ \sum_{k=0}^\infty 
     2^{-kp}\norm[T_1f][L^p(A_k)]^p \\
 &\les\ (1+|\lambda|)^{-2}\norm[V][]^p\sum_{k=0}^\infty 
         2^{-kp}\int_{D_k} |f(x)|^p\,dx
\end{aligned} $$
Interchange the summation and integral by Fubini's theorem.  At each point
$x \in \R^3$, $x\in D_k$ only if $k > \log(|x|/(2\lambda^{1/p}))$, 
so only these terms of the sum will be nonzero.
The resulting sum over $k$ is a geometric series with ratio less than 
$\frac12$, which can be estimated by the largest term.  Thus
$$ \begin{aligned}
\norm[T_1f][L^{p,-1}]^p\ &\les \ (1+|\lambda|)^{-2}\norm[V][]^p 
   \int_{\R^3} |f(x)|^p \min(2^p\lambda |x|^{-p},1)\,dx \\
&\les 2^p(1+|\lambda|)^{-1}\norm[V][]^p \int_{\R^3} |f(x)|^p(1+|x|)^{-p}\,dx
\end{aligned} $$
Taking $p^{\rm th}$ roots yields the desired conclusion.
\end{proof}
\begin{cor}
Suppose $V \in L^{\frac32(1+\eps)}(\R^3)\cap L^1(\R^3)$, and assume
$1 \le p \le 1+\eps$.  Then 
$VR_0^\pm(\lambda^2)$ is a bounded operator on $L^{p,1}(\R^3)$ whose norm
is controlled by $\eps$ and $\norm[V][]$ alone.
\end{cor}
\begin{proof}
This is the dual statement of Lemma~\ref{lem:weightedRV}, since V is 
real-valued and $L^{p',-\sigma}(\R^3)$ is the space dual to 
$L^{p,\sigma}(\R^3)$.
\end{proof}
\begin{cor} \label{cor:lambdadecay}
Suppose $V \in L^{\frac32(1+\eps)}(\R^3) \cap L^1(\R^3)$, and let 
$p \ge \frac{1+\eps}{\eps}$.  Then for all $\lambda\in\R$,
\begin{equation} \label{eq:lambdadecay}
\norm[(R_0^\pm(\lambda^2)V)^2f][L^{p,-1}] \les
  (1+|\lambda|)^{-1/p}\norm[V][]^2\norm[f][L^{p,-1}] 
\end{equation}
The dual operators satisfy the bound
\begin{equation}
\norm[(VR_0^\pm(\lambda^2))^2f][L^{p,1}] \les
  (1+|\lambda|)^{-1/p'}\norm[V][]^2\norm[f][L^{p,1}]\quad {\text for\ all\ 
   exponents}\ 1\le p\le 1+\eps \tag{\ref{eq:lambdadecay}'}
\end{equation}
\end{cor}
\begin{proof}
This is an estimate on $(T_1 + T_2)(T_1 + T_2)f$.  Any product which includes
$T_1$ will have the desired decay (or better) by Lemma~\ref{lem:weightedRV}.
On the other hand, 
$$ \begin{aligned} 
  \norm[T_2T_2f][L^\infty]\ \les \ \norm[V][] 
   \norm[T_2f][L^{\infty,-1}(\R^3\setminus B(0,\lambda^{1/p}))]
  \ &\les \ (1+|\lambda|)^{-1/p}\norm[V][]\norm[T_2f][L^\infty]  \\
    &\les \ (1+|\lambda|)^{-1/p}\norm[V][]^2\norm[f][L^{p,-1}]
\end{aligned} $$
\end{proof}

This is a crucial estimate for two reasons.  First, it guarantees convergence
of the Neumann series for $(I + R_0^\pm(\lambda^2)V)^{-1}$ for sufficiently
large $\lambda$, along with the uniform size bound
$$\limsup_{\lambda\to\infty} \norm[(I+R_0^\pm(\lambda^2)V)^{-1}][] \ \les \ 
  1 + \limsup_{\lambda\to\infty}\norm[R_0^\pm(\lambda^2)V][]$$
as measured in the operator norm on $L^{p,-1}(\R^3)$.  Second, we will
eventually perform an integration by parts in the $\lambda$ variable, whose
boundary terms will vanish because of \eqref{eq:lambdadecay}.

\section{Estimates on the perturbed resolvent}

Recall that the perturbed resolvent $R_V(z)$ is related to the $R_0(z)$ 
by the identity
\begin{equation}
R_V(z) = (I + R_0(z)V)^{-1}R_0(z) \tag{\ref{eq:resident}} 
\end{equation}
In order to prove that $R_V^\pm(\lambda^2)$ satisfies the same mapping
estimates as $R_0^\pm(\lambda^2)$, it therefore suffices to show that
$(I + R_0^\pm(\lambda^2)V)^{-1}$ is a bounded operator on the appropriate 
space.  As mentioned above, for large $\lambda$ this can be done easily
by expressing the inverse as a (convergent) power series.

If $\lambda$ is not large, invertability of $I+R_0^\pm(\lambda^2)V$ is 
established by a Fredholm-alternative argument.  One needs to verify two 
things: that $I + R_0^\pm(\lambda^2)V$ is a compact perturbation of the
identity, and that its null space contains no nonzero elements. 
This step will require the assumption that zero energy is neither an
eigenvalue nor a resonance, so we must first state a precise definition.
\begin{defi}
\label{def:res}
We say that a resonance occurs at zero if the equation $(I +R_0(0)V)f=0$ 
admits a distributional solution $f$ such that
$f\in L^{2,\sigma}(\R^3)\setminus L^2(\R^3)$ for every $\sigma < -\frac12$.
\end{defi}

\begin{lemma} \label{lem:compactness}
Suppose $V\in L^{\frac32(1+\eps)}(\R^3) \cap L^1(\R^3)$ and let 
$\frac{1+\eps}{\eps} \le p \le \infty$.  For any fixed $\lambda\in\R$,
$R_0^\pm(\lambda^2)V$ is a compact operator mapping $L^{p,-1}(\R^3)$ to itself.
By duality, $VR_0^\mp(\lambda^2)$ is a compact operator on $L^{p',1}(\R^3)$.
\end{lemma}
\begin{proof}
First consider the case where $V$ is bounded with maximum size $M$ 
and supported in the ball $B(0,R)$.  On the support of $V$, $f$ is 
integrable with bound
$\norm[f][L^1(\supp(V))] \les R^{1+3/p'}\norm[f][L^{p,-1}]$.
Then for all $|x| > 2R$,
$$|R_0^\pm(\lambda^2)Vf(x)|\ \les \ \big(|Vf| * \frac1{|\cdot|}\big)(x)
                          \ \les \ MR^{1+3/p'}\norm[f][L^{p,-1}]\,|x|^{-1}$$
Let $\psi$ be a smooth bump function with support in $B(0,2)$ 
so that $\psi(x) = 1$ whenever $|x| \le 1$, and define 
$\psi_{\tilde{R}}(x) = \psi(x/\tilde{R})$.
If $\tilde{R} > 2R$, a simple integration yields
$$\lim_{\tilde{R}\to\infty}
    \norm[(1-\psi_{\tilde{R}})R_0^\pm(\lambda^2)Vf][L^{p,-1}] \les 
  \lim_{\tilde{R}\to\infty}\big(MR^{1+3/p'}\norm[f][L^{p,-1}]\big)
     \tilde{R}^{1-3/p'} = 0$$

The resolvent tends to increase regularity; for Schwartz functions $f$ we have
$$ \begin{aligned}
(-\Lap + 1)R_0^\pm(\lambda^2)Vf &= (-\Lap-\lambda^2)R_0^\pm(\lambda^2)Vf
  + (1+\lambda^2)R_0^\pm(\lambda^2)Vf \\
 &= Vf +  (1+\lambda^2)R_0^\pm(\lambda^2)Vf
\end{aligned} $$
which implies that $\norm[(-\Lap+1)R_0^\pm(\lambda^2)Vf][L^{p,-1}] \les
  \norm[f][L^{p,-1}]$.  Boundedness of $V$ is also used in this step.
Taking limits, the inequality can be extended to all $f\in L^{p,-1}$.  

On the compact set $\{|x| \le 2\tilde{R}\}$ the norms $L^p$ and $L^{p,-1}$
are equivalent. 
Therefore $\psi_{\tilde{R}} R_0^\pm(\lambda^2)V$ is a continuous map
from $L^{p,-1}$ to the Sobolev space $W^{2,p}(B(0,2\tilde{R}))$, which embeds
compactly into $L^p(B(0,2\tilde{R}))$, and hence also $L^{p,-1}(\R^3)$,
by Rellich's theorem.

We have shown that $R_0^\pm(\lambda^2)V$ is a norm-limit of the operators
$\psi_{\tilde{R}}R_0^\pm(\lambda^2)V$ as $\tilde{R} \to \infty$, and that
each element of this approximating sequence is compact.  The set of compact
linear operator is closed in the norm topology, so $R_0^\pm(\lambda^2)V$
must be compact as well.

For general potentials $V \in L^{\frac32(1+\eps)}(\R^3) \cap L^1(\R^3)$,
it is possible write $V$ as a norm-limit of bounded
functions $V_n$ with compact support. For each $n = 1,2,\ldots$,  
$R_0^\pm(\lambda^2)V_n$ is a compact operator.  The lemma is now proved by
another limiting argument, this time with the help of 
lemma~\ref{lem:weightedRV}.
\end{proof}

\begin{lemma} \label{lem:fredholm}
Let $V \in L^{\frac32(1+\eps)}(\R^3) \cap L^1(\R^3)$, and $1\le p\le 1+\eps$.
Assume that zero is neither an eigenvalue nor a resonance of $(-\Lap + V)$.
Then $(I + VR_0^\pm(\lambda^2))^{-1}$ exists as a bounded linear map on 
$L^{p,1}(\R^3)$ for all $\lambda \in \R$.  By duality, 
$(I + R_0^\pm(\lambda^2)V)^{-1}$ 
exists as a bounded operator on $L^{p,-1}(\R^3)$.
\end{lemma}
\begin{proof}
By lemma \ref{lem:compactness} and the Fredholm alternative, 
$I + VR_0^{\pm}(\lambda^2)$ will fail be invertible only if there
exists a function $g\in L^{p,1}(\R^3)$ satisfying $g = -VR_0^\pm(\lambda^2)g$.
In fact any such solution $g$ must possess greater regularity than the assumed
$g \in L^p_{\rm loc}$.  This is seen by iterating the map $VR_0^\pm(\lambda^2)$.  

First note that any function in $L^{q,1}(\R^3), q < \frac32$, is integrable.
Decompose $R_0^\pm(\lambda^2) = S_1 + S_2$
in the following manner:  For $x \in A_k,\ k=1,2,\ldots$, let 
$$ \begin{aligned} S_1f(x) &= R_0^\pm(\lambda^2)\chi_{\{|x| > 2^{k-2}\}}f(x) \\
                   S_2f(x) &= R_0^\pm(\lambda^2)\chi_{\{|x| \le 2^{k-2}\}}f(x)
\end{aligned} $$
and $S_1g(x) = R_0^\pm(\lambda^2)g(x)$ if $x \in A_0$.  Here, as in
lemma~\ref{lem:weightedRV}, $A_k$ denotes an annulus where $|x| \sim 2^k$.
One immediately obtains a pointwise estimate for $S_2$, namely
$|S_2f(x)| \les \norm[f][L^1] (1+|x|)^{-1}$.

We will see that $S_1$ is a bounded map from $L^{q,1}(\R^3)$ to 
$L^{r,1}(\R^3)$, where $\frac1r = \frac1q - \frac23$ is the exponent given
by fractional integration.  The calculation is similar to the one in
lemma~\ref{lem:weightedRV}, with one additional step to deal with the fact
that $r \not= q$.

\begin{equation} \label{eq:weightedR} \begin{aligned}
\norm[S_1g][L^{r,1}]\ \sim \ \Big(\sum_{k=0}^\infty 2^{kr} 
 \norm[S_1g][L^r(A_k)]^r \Big)^{1/r} \ &\les \ 
 \Big(\sum_{k=0}^\infty 2^{kr} \big(\int_{|x| \ge 2^{k-2}} |g(x)|^q\,dx 
    \big)^{r/q} \Big)^{1/r}  \\
  &\le \Big(\int_{\R^3} |g(x)|^q \big(\sum_{k\le \log 4|x|} 2^{kr}\big)^{q/r}
      dx \Big)^{1/q} \\
  &\les \big(\int_{\R^3}|g(x)|^q |x|^q\, dx\big)^{1/q} \ = \ \norm[g][L^{q,1}]
\end{aligned}
\end{equation}
The exchange of summation and integration is done via Minkowski's inequality,
noting that $r > q$.

Putting the two pieces $S_1$ and $S_2$ together, we conclude that 
$R_0^\pm(\lambda^2)$ 
is a bounded map from $L^{q,1}(\R^3)$ to $L^{r,1}(\R^3) + L^{\infty,1}(\R^3)$.
Therefore, if $g \in L^{q,1}(\R^3)$, one can bootstrap in two directions:
\begin{equation} \label{eq:bootstrap}
VR_0^\pm(\lambda^2)g \in L^{1,1}(\R^3) \qquad{\rm and}\qquad
VR_0^\pm(\lambda^2)g \in L^{\tilde{q},1}(\R^3),\quad {\rm where}\ \ 
 \frac1{\tilde{q}} = \frac1q - \frac{2\eps}{3(1+\eps)}
\end{equation}
by estimating $\norm[V][]$ in $L^{r'}\cap L^1$ and 
$L^{\frac32(1+\eps)}\cap L^{\tilde{q}}$, respectively.

This process terminates once it is established that $g \in L^{1,1}(\R^3)
\cap L^{\frac32^+,1}(\R^3)$.  Consequently, $g \in L^{\frac32^+}(\R^3)\cap
L^{\frac32^-}(\R^3)$, and 
$R_0^\mp(\lambda^2)g \in L^{\infty}(\R^3)$, which embeds naturally in 
$L^{\infty,-1}(\R^3)$.  The pairing of functions in dual spaces
$$\la R_0^\pm(\lambda^2)g, g\ra \ = \ 
    -\la R_0^\pm(\lambda^2)g,VR_0^\pm(\lambda^2)g\ra$$
is then well-defined.  Furthermore, since $V$ is asumed to be real-valued,
the expression on the right side has no imaginary part.  On the other hand,
by Parseval's identity
\begin{equation} \label{eq:restriction}
\Im\la \R_0^\pm(\lambda^2)g, g\ra \ = \ 
  \lim_{\eps\to 0}\Im\la R_0(\lambda^2\pm i\eps)g, g\ra \ = \ 
  \pm C\lambda \int_{S^2} |\hat{g}(\lambda\omega)|^2 \sigma(d\omega)
\end{equation}
where $C \not= 0$ is a constant and $\sigma(d\omega)$ is surface measure
on the unit sphere in $\R^3$.  It follows that $\hat{g} = 0$ on $\lambda S^2$,
in the sense of $L^2$ functions.

One of the underlying principles in Agmon~\cite{agmon} is that the resolvent
$R_0^\pm(\lambda^2)$ has special mapping properties when applied to functions
whose Fourier transform vanishes on the sphere radius $\lambda$.  This in turn
leads to improved estimates on the decay of 
$g = -VR_0^\pm(\lambda^2)g$.  We quote one such statement from the literature:
\begin{lemma}[\cite{GS2}, section 4]
Let $f$ be a function in $L^1(\R^3)$ such that $\hat{f} = 0$ on the unit 
sphere.  Then
$$ \norm[R_0^\pm(1)f][L^2] \le \frac1{\sqrt{8\pi}} \norm[f][L^1] $$
\end{lemma}
\noindent
Scaling considerations dictate that if $\hat{g} = 0$ on $\lambda S^2$, then
$\norm[R_0^\pm(\lambda^2)g][L^2] \le (8\pi\lambda)^{-1/2} \norm[g][L^1]$.

Returning to the proof of lemma \ref{lem:fredholm}, if $\lambda \not= 0$ and 
$g \in L^{p',1}(\R^3)$ is a solution to $(I + VR_0^\pm(\lambda^2))g = 0$,
then $f = R_0^\pm(\lambda^2)g$ must be an $L^2$~eigenfunction of $-\Lap + V$.
The bootstrapping procedure for $g$ shows that $f \in W^{2,3/2}_{\rm loc}(\R^3)
\subset W^{1,2}_{\rm loc}(\R^3)$.  By assumption $V \in L^{\frac32}(\R^3)$, 
so all the hypotheses of Theorem~\ref{thm:IJ} are satisfied.
One concludes that $f = 0$, and $g = -Vf = 0$, as desired.

In the case $\lambda = 0$, the expression in \eqref{eq:restriction} is
trivially zero, so $\hat{g}$ does not satisfy 
any additional hypotheses.  The resolvent $R_0(0)$ is a bounded map from 
$L^1(\R^3)$ to $L^{2,-\sigma}(\R^3)$ for any $\sigma > \frac12$, however, 
so $f = R_0(0)g$ is a distributional solution of $-\Lap + V$ which lies in 
every space $L^{2,-\sigma}(\R^3), \sigma > \frac12$.  The assumption that 
zero energy is neither an eigenvalue nor a resonance requires that $f = 0$, 
thus $g = -Vf = 0$ as well.
\end{proof}

\noindent
{\bf Remark.}
To be precise, $R_0^\pm(\lambda^2)$ maps $L^{1,1}(\R^3)$ to 
weak-$L^{3,1}(\R^3)$.  The calculation in \eqref{eq:weightedR}, with the
appropriate cosmetic changes, is used to bound 
$\big({\rm meas}\{|R_0^\pm(\lambda^2)g(x)| > h/(1+|x|)\}\big)^{1/3}h$
uniformly for all $h > 0$.  

The bootstrapping estimates in \eqref{eq:bootstrap} will be needed in 
Section~\ref{sec:calculations} in the following form:
\begin{prop} \label{prop:weighted(VR)^k}
Let $V \in L^{\frac32(1+\eps)}(\R^3) \cap L^1(\R^3)$.  Then 
\begin{align}  
\norm[VR_0^\pm(\lambda^2)f][L^{2/(2-\eps),1}] &\les\norm[f][L^{1,1}] \notag\\
\norm[VR_0^\pm(\lambda^2)f][L^{1,1}] &\les \norm[f][L^{2/(2-\eps),1}] \notag\\
  \norm[(VR_0^\pm(\lambda^2))^{k+3}f][L^{1,1}] &\les 
(1+|\lambda|)^{-k\eps/4} \norm[f][L^{1,1}] \label{eq:L11powerdecay}
\end{align}
\end{prop}
\begin{proof}
The first two inequalties are precisely what is stated in \eqref{eq:bootstrap}.
The last line combines these with (\ref{eq:lambdadecay}').
Neither the choice of exponent $p = \frac2{2-\eps}$ nor the power of
decay in $\lambda$ are intended to be sharp.  For our purposes it will
only matter that $k$ can be chosen large enough to make $k\eps/4 > 1$.
\end{proof}

\begin{prop} \label{prop:continuousRV}
Let $V \in L^{\frac32(1+\eps)}(\R^3) \cap L^1(\R^3)$.
The family of linear maps $VR_0^\pm(\lambda^2)$, considered with respect to
the operator norm on $L^{p,1}(\R^3),\ 1\le p<\frac32$, depends continuously 
on the parameter $\lambda$.  By duality, $R_0^\mp(\lambda^2)V$ is a continuous
family of maps on $L^{p,-1},\ p>3$.
\end{prop}
\begin{proof}
Any function $f\in L^{p,1}(\R^3)$ is also integrable, with $\norm[f][L^1]
\les \norm[f][L^{p,1}]$.  It is then possible to differentiate under the
integral sign to obtain
\begin{equation*}
\Big|\dfrac{d}{d\lambda}R_0^\pm(\lambda^2)f(x)\Big| \ = \ 
\Big|\int_{\R^3} (\mp4\pi i)^{-1} e^{\pm i \lambda |x-y|}f(y)\,dy\Big| \ \les 
\ \norm[f][L^{p,1}]
\end{equation*}
If $V$ is bounded and has compact support, then $V \in L^{p,1}$ and
$\frac{d}{d\lambda}\big[VR_0^\pm(\lambda^2)\big]$ will be a bounded
operator on $L^{p,1}(\R^3)$ uniformly in $\lambda$, which implies continuity.

For general potentials $V$, approximate $V$ by a bounded, compactly supported
potential $V'$ so that $\norm[V-V'][] < \eps$.  Then 
$\norm[(V-V')R_0^\pm(\lambda^2)f][L^{p,1}] < C\eps\norm[f][L^{p,1}]$
where $C<\infty$ is the constant in lemma~\ref{lem:weightedRV}.  
Assume that $V$ is supported in the ball $B(0,R)$ and that $\sup_x |V(x)| = M$.
For every $|\nu|\in \R$,
the operator $R_0^\pm((\lambda+\nu)^2) - R_0^\pm(\lambda^2)$
is bounded from $L^1(\R^3)$ to $L^\infty(\R^3)$ with norm $\frac{|\nu|}{4\pi}$,
therefore
$$\norm[V'{[R_0^\pm((\lambda+\nu)^2) - R_0^\pm(\lambda^2)]}f][L^{p,1}] \les 
  MR^4|\nu|\norm[f][L^{p,1}]$$
By the triangle inequality
$$\liminf_{\nu\to 0}\norm[V{[R_0^\pm((\lambda+\nu)^2)-R_0^\pm(\lambda^2)]}f][
L^{p,1}] < 2C\eps\norm[f][L^{p,1}]$$
\end{proof}
\begin{lemma} \label{lem:supInverse}
Let $V \in L^{\frac32(1+\eps)}(\R^3) \cap L^1(\R^3)$ and assume that zero
energy is neither an eigenvalue nor a resonance.  Then
\begin{equation}
\sup_{\lambda\in\R} \norm[(I+VR_0^\pm(\lambda^2))^{-1}][L^{p,1}\to L^{p,1}]
  < \infty \qquad {\it for\ all}\ \ 1 \le p \le 1+\eps
\end{equation}
\end{lemma}
\begin{proof}
Consider the case $p > 1$.
By corollary~\ref{cor:lambdadecay}, there exists $\lambda_0<\infty$ 
so that the operator norm of 
$(VR_0^\pm(\lambda^2))^2$ will be less than $\frac12$ for all 
$|\lambda| > \lambda_0$.
For these large values of $\lambda$, the Neumann series
$$(I + VR_0^\pm(\lambda^2))^{-1} = \sum_{k=0}^\infty (VR_0^\pm(\lambda^2))^{2k}
  (I + VR_0^\pm(\lambda^2))$$
converges geometrically and has norm controlled by 
$(1 + \norm[VR_0^\pm(\lambda^2)][]) \les 1 + \norm[V][]$.
At every point $\lambda \in \R$, 
lemma~\ref{lem:fredholm} and proposition~\ref{prop:continuousRV} 
and the continuity of inverses guarantee that $(I + VR_0^\pm(\lambda^2))^{-1}$
is norm-continuous in $\lambda$.  Thus it is bounded on the compact set 
$[-\lambda_0, \lambda_0]$.

In the case $p=1$, we claim that 
$\norm[(VR_0^\pm(\lambda^2))^2][L^{1,1}\to L^{1,1}]$ vanishes as 
$|\lambda| \to \infty$.
A substantially similar result appears in
\cite{DanPie}, which we reproduce below with the necessary modifications.
Assuming this fact, the remaining steps of the above argument follow
immediately.
\end{proof}
\begin{proposition}
Suppose $V \in L^{\frac32(1+\eps)}(\R^3)\cap L^1(\R^3)$.  Then
\begin{equation} \label{eq:L11decay}
\lim_{\lambda\to\infty} 
\norm[(VR_0^\pm(\lambda^2))^2][L^{1,1}\to L^{1,1}] = 0
\end{equation}
\end{proposition}
\begin{proof}
Suppose $V$ is supported in the ball $B(0,R)$ and satisfies $|V(x)| < M$.
Then 
$\norm[VR_0^\pm(\lambda^2)f][L^{4/3}] \les 
R^{\frac54}\norm[VR_0^\pm(\lambda^2)f][L^3_{\rm weak}] \les 
MR^{\frac54}\norm[f][L^{1,1}]$.  
It follows from proposition~\ref{prop:limabs} that
$$\norm[(VR_0^\pm(\lambda^2))^2f][L^{1,1}] \les 
R^{\frac{13}4}\norm[(VR_0^\pm(\lambda^2))^2f][L^4] \les 
M^2R^{\frac92}\lambda^{-1/2}\norm[f][L^{1,1}]$$
Stronger decay estimates are possible, but we are not interested here in 
optimality.

For general potentials $V$, write $V = V_1 + V_2$ with $V_1$ bounded
and compactly supported and $\norm[V_2][] < \epsilon$.  Then
$$\norm[(VR_0^\pm(\lambda^2))^2][] \le \norm[(V_1R_0^\pm(\lambda^2))^2][]
+ \norm[V_1R_0^\pm(\lambda^2)V_2R_0^\pm(\lambda^2)][] +
 \norm[V_2R_0^\pm(\lambda^2)VR_0^\pm(\lambda^2)][]$$
All three terms on the right-hand side are smaller than $\eps$ when
$\lambda$ is sufficiently large.
\end{proof}

{\bf Remark.}  It is also true that the operators 
$(I + VR_0^\pm(\lambda^2))^{-1}$ are uniformly bounded on the unweighted spaces
$L^p(\R^3), 1\le p < \frac32$.  The proof follows the same 
Fredholm-alternative argument, but uses \eqref{eq:freeVR} in place of
\eqref{eq:lambdadecay} and similar substitutions.  The details (with more
restrictive hypotheses on $V$) for $p=1$ can be found in \cite{DanPie}
and for $p = \frac43$ in \cite{GS2}.

The primary condition on $V$, that $V \in L^{\frac32(1+\eps)}(\R^3)\cap
L^1(\R^3)$, is translation-invariant.  Indeed, if $V_y(x) = V(x-y)$ is
any translate of $V$, then $\norm[V_y][] = \norm[V][]$.  The second condition,
that zero energy is neither an eigenvalue nor resonance for $-\Lap + V$,
is also preserved under translation.  The norm of functions in $L^{p,1}(\R^3)$,
however, is clearly affected by translations and cannot even be bounded
uniformly.  Nevertheless a translation-invariant statement of 
lemma~\ref{lem:supInverse} is still possible.
\begin{lemma} \label{lem:SupInverse}
Let $V \in L^{\frac32(1+\eps)}(\R^3) \cap L^1(\R^3)$ and assume that zero
energy is neither an eigenvalue nor a resonance.  Then
\begin{equation}
\sup_{y\in\R^3} \sup_{\lambda\in\R} \norm[(I+V_yR_0^\pm(\lambda^2))^{-1}][
L^{p,1}\to L^{p,1}]
  < \infty \qquad {\it for\ all}\ \ 1 \le p \le 1+\eps
\end{equation}
\end{lemma}
\begin{proof}
The mapping $y \mapsto V_y$ is uniformly continuous, that is
$\norm[V_y - V_z][] < \delta$ for every pair of points with $|y-z| > \delta'$.
Consequently, the family of operators $V_yR_0^\pm(\lambda^2)$ are uniformly
continuous with respect to variation in the $y$ parameter.  Meanwhile, by
proposition~\ref{prop:continuousRV} this family is also continuous with respect
to variation in the $\lambda$ parameter.  It follows that
$$(y,\lambda) \in \R^3\times\R \mapsto V_yR_0^\pm(\lambda^2) \in 
{\mathcal B}(L^{p,1}(\R^3))$$
is continuous on $\R^3\times\R$.

The decay estimates \eqref{eq:lambdadecay} and \eqref{eq:L11decay} hold
uniformly over all translations of $V$.  Thus there exists $\lambda_0 <\infty$
such that $\norm[(V_yR_0^\pm(\lambda^2))^2][] < \frac12$ for all 
$|\lambda| > \lambda_0$ and all $ y\in\R^3$.  As in the proof of 
lemma~\ref{lem:supInverse}, the operator norm of 
$(I + V_yR_0^\pm(\lambda^2))^{-1}$ is controlled uniformly by 
$1 + \norm[V][]$ at these points.

Suppose $V$ is supported in the ball $B(0,r)$ and $|y| > 3r$.  Given a 
function $f \in L^{p,1}(\R^3)$, let $f_1 = \chi_{B(y,2r)}f$ and
$f_2 = f - f_1$.  By construction, $f + V_yR_0^\pm(\lambda^2)f = f_2$
outside the ball $B(y,2r)$, thus $\norm[f + V_yR_0^\pm(\lambda^2)f][L^{p,1}]
\ge \norm[f_2][L^{p,1}]$.

Within $B(y,2r)$, we have that $f + V_yR_0^\pm(\lambda^2)f = f_1 
+ V_yR_0^\pm(\lambda^2)(f_1 + f_2)$.  Thus
$$\norm[f+V_yR_0^\pm(\lambda^2)f][L^{p,1}] \ge 
\norm[f_1 + V_yR_0^\pm(\lambda^2)f_1][L^{p,1}] 
  - \norm[V_yR_0^\pm(\lambda^2)f_2][L^{p,1}]$$
Since every point $x \in B(y,2r)$ satisfies $r < |x| < 5r$, the weighted and
unweighted norms are equivalent.  
There exists a constant $A > 0$ such that
$\norm[g + VR_0^\pm(\lambda^2)g][L^p] \ge A\norm[g][L^p]$ for all $\lambda
\in \R$ and every $g \in L^p(\R^3)$.  This is equivalent to the uniform
boundedness of $(I + VR_0^\pm(\lambda^2))^{-1}$ as operators on $L^p(\R^3)$.
By translation invariance, the same
estimate holds if $V$ is replaced by any $V_y$.  It follows that
$$\norm[f_1 + V_yR_0^\pm(\lambda^2)f_1][L^{p,1}] \ge 
  \frac{A}5\norm[f_1][L^{p,1}]$$
since the functions on both sides are supported in $B(y,2r)$.
For $f_2$, the crude estimate $|R_0^\pm(\lambda^2)f_2(x)| \le (4\pi r)^{-1}
\norm[f_2][L^1] \les (4\pi r)^{-1}\norm[f_2][L^{p,1}]$ is valid at all
$x\in B(y,2r)$.  This suffices to show that
$$\norm[V_yR_0^\pm(\lambda^2)f_2][L^{p,1}] \le C\norm[V_y][L^p]
  \norm[f_2][L^{p,1}]$$
Applying the triangle inequality to $f = f_1 + f_2$, we conclude that
$$\norm[f + V_yR_0^\pm(\lambda^2)f][L^{p,1}] \ge \max\big(\norm[f_2][],\ 
\frac{A}5\norm[f][] - (\frac{A}5 + C\norm[V][])\norm[f_2][]\big) \ge 
\frac{A}{A + 5C\norm[V][] + 5}\norm[f][L^{p,1}]$$

Here we are taking advantage of the fact that $\norm[V][] = \norm[V_y][]$ 
Observe that none of the constants in this inequality depend on the 
size or support of $V$.  Given an arbitrary potential 
$V \in L^{\frac32(1+\eps)}\cap L^1$,
it is then possible to choose $V_r = \chi_{|x|<r}V$ so that 
$\norm[V-V_r][] \les \frac{A}{2(A + 5C\norm[V][] + 5)}$.  By 
lemma~\ref{lem:weightedRV} and the triangle inequality,
$$\norm[f + V_yR_0^\pm(\lambda^2)f][L^{p,1}] \ge
  \frac{A}{2(A + 5C\norm[V][] + 5)}\norm[f][L^{p,1}]$$
for all $|y| > 3r$.

Having established a uniform bound on $(I + V_yR_0^\pm(\lambda^2))^{-1}$
for all $|\lambda| > \lambda_0$ and for all $|y| > 3r$, only a compact region
of $\R^3\times\R$ remains to be considered.  However 
$(I + V_yR_0^\pm(\lambda^2))^{-1}$ is a continuous function of $(y,\lambda)$,
hence it is bounded on this domain as well.
\end{proof}

\section{Calculations} \label{sec:calculations}

Our goal at this point is to prove the estimate
\begin{equation} \begin{aligned}
\int_0^\infty e^{it\lambda^2}\lambda \Big\la 
  \big[R_V^+(\lambda^2)(VR_0^+(\lambda^2))^{m+2}
 -\ R_V^-(\lambda^2)(VR_0^-(\lambda^2))^{m+2} \big]f, g \Big\ra\, d\lambda \\
=\ \int_0^\infty e^{it\lambda^2}\lambda \la A(\lambda)f, g\ra\,d\lambda
\ \les\ |t|^{-\frac32}\norm[f][1]\norm[g][1]
\end{aligned}
\end{equation}
This will be true if and only if the operator
$$\int_0^\infty e^{it\lambda^2}\lambda\, A(\lambda)\,d\lambda$$
is a well defined map from $L^1$ to $L^\infty$ whose operator norm is
controlled by $|t|^{-3/2}$.

\begin{lemma} \label{lem:Adecay}
For sufficiently large values of m, $\lim_{\lambda\to\infty}
\norm[A(\lambda)][] = 0$ as a map from $L^1(\R^3)$ to $L^\infty(\R^3)$.
\end{lemma}
\begin{proof}
Decompose $A(\lambda)$ into a telescoping series
$$\begin{aligned}
A(\lambda)\ =\ \Big[(R_V^+(\lambda^2) - R_V^-(\lambda^2))
  &(VR_0^+(\lambda^2))^{m+2} + \\ 
 &(R_V^-(\lambda^2)V)\sum_{k=0}^{m+1} (R_0^-(\lambda^2)V)^k
  (R_0^+(\lambda^2) - R_0^-(\lambda^2))(VR_0^+(\lambda^2))^{m+1-k}\Big]
\end{aligned} $$
The difference $R_0^+(\lambda^2) - R_0^-(\lambda^2)$ is precisely a convolution
with the kernel $\frac{i\,sin(\lambda |x|)}{2\pi |x|}$.  This maps $L^1$ to
$L^\infty$ with operator norm proportional to $\lambda$.  The difference of
perturbed resolvents has a similar bound, by the identity
\begin{equation}
R_V^+(\lambda^2) - R_V^-(\lambda^2) = (I+R_0^+(\lambda^2)V)^{-1}
   (R_0^+(\lambda^2) - R_0^-(\lambda^2))(I + VR_0^-(\lambda^2))^{-1}
\end{equation}
Choose $m$ so that $(m-3)\eps > 1$.  By \eqref{eq:(RV)^k}, 
each of the $m+2$ terms is a bounded map from $L^1$ to $L^\infty$ with norm
controlled by $\lambda(1+|\lambda|)^{-(m-3)\eps}$.  
\end{proof}

We may then integrate by parts to obtain
\begin{equation}\label{eq:IBP}
\int_0^\infty e^{it\lambda^2}\lambda\, A(\lambda)\,d\lambda
= -\frac{1}{2it}\int_0^\infty e^{it\lambda^2} A'(\lambda)\,d\lambda
\end{equation}
The boundary term at infinity vanishes by lemma~\ref{lem:Adecay}.
The boundary term at $\lambda=0$ vanishes because $R_0^+(0) = R_0^-(0)$.
From this point forward, cancellation involving 
$R_0^+(\lambda^2) - R_0^-(\lambda^2)$ will not play a major role; this allows
us to express $A'(\lambda)$ in a less cumbersome manner.  Recall that 
$$A'(\lambda) = \dfrac{d}{d\lambda}\big[R_V^+(\lambda^2)
(VR_0^+(\lambda^2))^{m+2}\big] -\ \dfrac{d}{d\lambda}\big[R_V^-(\lambda^2)
(VR_0^-(\lambda^2))^{m+2}\big] \ := \ B^+(\lambda) - B^-(\lambda)$$

For all $\lambda > 0$, the resolvent $R_V^+(\lambda^2)$ may be defined via 
two different limits, since $R_V(\lambda^2+i0)\ =\ R_V((\lambda+i0)^2)$.  
If zero energy is neither a resonance nor an eigenvalue,
the latter expression admits an analytic continuation into the half-plane 
$\Im \lambda > 0$, with continuous extension to the boundary satisfying
$$R_V((-\lambda+i0)^2) = R_V((\lambda-i0)^2)$$
The same expression is of course true for the free resolvent as well.
A similar analytic extension exists for $B^+(\lambda)$, with the boundary
identity
$B^+(-\lambda) = -B^-(\lambda)$.
The change in sign is a result of differentiation with respect to 
$\lambda$.  The right-hand integral in \eqref{eq:IBP} may now be rewritten as
$$
-\frac1{2it}\int_{-\infty}^\infty e^{it\lambda^2} B^+(\lambda)\, d\lambda
$$

The proof of theorem~\ref{thm:dispersive} concludes with an
estimate for this integral.

\begin{lemma}
With $V \in L^{\frac32(1+\eps)}(\R^3) \cap L^1(\R^3)$ and $B^+(\lambda)$
defined as above,
$$\Big\| \int_{-\infty}^\infty e^{it\lambda^2}B^+(\lambda) f\Big\|_{L^\infty}
\les |t|^{-\frac12} \norm[f][L^1] $$
for all functions $f\in L^1(\R^3)$.
\end{lemma}
\begin{proof}
We may express $B^+(\lambda)$ as an integral operator whose kernel is given
formally by the expression
\begin{equation} \label{eq:Bplus}
B^+(\lambda,x,y) = \dfrac{d}{d\lambda}\Big\la(I+R_0^+(\lambda^2)V)^{-1}
 R_0^+(\lambda^2)(VR_0^+(\lambda^2))^m V(\cdot)\frac{e^{i\lambda|\cdot - x|}}
{4\pi|\cdot-x|}, V(\cdot)\frac{e^{-i\lambda|\cdot - y|}}{4\pi|\cdot - y|}
\Big\ra
\end{equation}
The inner product as written above may not be well-defined because of the
local singularites in $V$.  If, however, the derivative 
is brought inside and applied according to the Leibniz
rule, then each term will be finite.  Essentially this is because
$\frac{d}{d\lambda}R_0^+(\lambda^2)$ is a uniformly bounded operator
from $L^1(\R^3)$ to $L^{\infty}(\R^3)$, leading to a pairing between
the dual spaces $L^1(\R^3)$ and $L^\infty(\R^3)$.

It is sufficient to show that $B^+(\lambda,x,y)$ is uniformly bounded, and
$\int_\R \big|\frac{d}{d\lambda}\big[e^{-i\lambda|y-x|}B^+(\lambda,x,y)\big]
\big|\,d\lambda$ is bounded uniformly in $x$ and $y$.  Then the lemma follows
from a stationary-phase argument, estimating the size of the integral
$$\int_{-\infty}^\infty e^{it(\lambda + \frac{|y-x|}{2t})^2}
   \big[e^{-i\lambda|y-x|}B^+(\lambda,x,y)\big]\,d\lambda$$

The derivative in \eqref{eq:Bplus} can fall in any of $m+4$ locations,
leading to a sum of four terms:
\begin{align*}
16\pi^2B^+(&\lambda,x,y) = \\
&i\Big\la(I+R_0^+(\lambda^2)V)^{-1}(R_0^+(\lambda^2)V)^{m+1}
e^{i\lambda|\cdot-x|},\ \frac{V(\cdot)e^{-i\lambda|\cdot-y|}}{|\cdot-y|}
\Big\ra   \tag{\ref{eq:Bplus}a} 
\\
+\ &\sum_{k=0}^m \Big\la(I+R_0^+(\lambda^2)V)^{-1}(R_0^+(\lambda^2)V)^k
\Big[\dfrac{d}{d\lambda}R_0^+(\lambda^2)\Big]
(VR_0^+(\lambda^2))^{m-k}\frac{V(\cdot)e^{i\lambda|\cdot-x|}}{|\cdot-x|},\ 
\frac{V(\cdot)e^{-i\lambda|\cdot-y|}}{|\cdot-y|}\Big\ra \tag{\ref{eq:Bplus}b}
\\
-\ &\Big\la(I+R_0^+(\lambda^2)V)^{-1}\Big[\dfrac{d}{d\lambda}R_0^+(\lambda^2)
\Big](I+VR_0^+(\lambda^2))^{-1}(VR_0^+(\lambda^2))^{m+1}
\frac{V(\cdot)e^{i\lambda|\cdot-x|}}{|\cdot-x|},\ 
\frac{V(\cdot)e^{-i\lambda|\cdot-y|}}{|\cdot-y|}\Big\ra \tag{\ref{eq:Bplus}c}
\\
+\ &i\Big\la(I+VR_0^+(\lambda^2))^{-1}(VR_0^+(\lambda^2))^{m+1}
\frac{V(\cdot)e^{i\lambda|\cdot-x|}}{|\cdot-x|},\ e^{-i\lambda|\cdot-y|}
\Big\ra \tag{\ref{eq:Bplus}d}
\end{align*}
The formula in (\ref{eq:Bplus}c) is a consequence of the chain rule
$\frac{d}{d\lambda}M^{-1}(\lambda) = 
M^{-1}(\lambda)\big[\frac{d}{d\lambda}M(\lambda)\big]M^{-1}(\lambda)$
for operator-valued functions,
and also the commutator relation $V(I+R_0^+(\lambda^2)V)^{-1}
= (I+VR_0^+(\lambda)^2)^{-1}V$.

Each of the four terms is bounded uniformly in $(\lambda,x,y)$ by some 
combination of \eqref{eq:freeVR}, \eqref{eq:freeRV}, H\"older's inequality,
and the following observations:  

$\sup_x\norm[e^{i\lambda|\cdot-x|}][\infty] = 1$.  

$\sup_x\norm[V(\cdot)|\cdot-x|^{-1}][1] \les \norm[V][]$.

$\sup_\lambda\norm[(I+VR_0^+(\lambda^2))^{-1}][1\to 1] < \infty$.  The proof
is essentially identical to that of lemma~\ref{lem:supInverse}.

$\sup_\lambda \big\|\frac{d}{d\lambda}R_0^+(\lambda^2)\big\|_{1\to\infty}
= (4\pi)^{-1}$.

\smallskip \noindent
We now turn our attention to the second assertion, that
$\sup_{(x,y)}\int_\R \big|\frac{d}{d\lambda}\big[e^{-i\lambda|y-x|}
B^+(\lambda,x,y)\big]\big|\,d\lambda < \infty$.  In fact we will prove
the pointwise estimate
\begin{equation} \label{eq:pointwise}
\big|\frac{d}{d\lambda}\big[e^{-i\lambda|y-x|}B^+(\lambda,x,y)\big]\big|
\les (1+|\lambda|)^{-(m-6)\eps/4} \quad {\it for\ all}\ m \ge 8.
\end{equation}
so that it suffices to choose $m > \frac4\eps + 6$.
For the sake of brevity, we will only calculate explicitly the derivatives
associated to a typical term in the expression (\ref{eq:Bplus}b).
The same techniques apply equally well to each of the other terms.

Suppose the derivative falls anywhere except on the already-differentiated
resolvent, a typical example being
\begin{multline*}
\Big\la(I+R_0^+(\lambda^2)V)^{-1}(R_0^+(\lambda^2)V)^\ell
\Big[\dfrac{d}{d\lambda}R_0^+(\lambda^2)\Big]V(R_0^+(\lambda^2)V)^{k-\ell-1} \\
\Big[e^{-i\lambda|x-y|}\dfrac{d}{d\lambda}R_0^+(\lambda^2)\Big]
(VR_0^+(\lambda^2))^{m-k}\frac{V(\cdot)e^{i\lambda|\cdot-x|}}{|\cdot-x|},\
\frac{V(\cdot)e^{-i\lambda|\cdot-y|}}{|\cdot-y|}\Big\ra
\end{multline*}
Using \eqref{eq:(RV)^k} and the four observations listed above, this term
is seen to be less than $(1+|\lambda|)^{-(m-7)\eps}$.
Of particular note here is the fact that
multiplication by $V$ is a bounded map between $L^\infty(\R^3)$ and
$L^1(\R^3)$.
 
The case where the derivative falls on $(I + R_0^+(\lambda^2)V)^{-1}$ has only
superficial differences, since the operator
$$\dfrac{d}{d\lambda}(I+R_0^+(\lambda^2)V)^{-1} = (I + R_0^+(\lambda^2)V)^{-1}
\Big[\dfrac{d}{d\lambda} R_0^+(\lambda^2)\Big] (I + VR_0^+(\lambda^2))^{-1}V$$ 
is still bounded on $L^\infty(\R^3)$ uniformly in $\lambda$.

In order to address the case where both derivatives fall in the same
place, we use estimates in $L^{p,\sigma}_x(\R^3)$, the weighted norm space
defined by
$$\norm[f][L^{p,\sigma}_x] := \norm[(1+|\cdot - x|)^\sigma f][L^p]$$
This cannot easily be avoided, as the kernel of $\dfrac{d^2}{d\lambda^2}
R_0^+(\lambda^2)$ experiences polynomial growth in the spatial variables.
Note that $L^{p',-\sigma}_x(\R^3)$ is the dual space to $L^{p,\sigma}_x(\R^3)$ 
for any $1 \le p < \infty, \sigma \in \R$.

It is clear by translation that the action of $VR_0^+(\lambda^2)$ on
$L^{p,1}_x$ is equivalent to that of $V_{-x}R_0^+(\lambda^2)$ acting on 
$L^{p,1}$.  The bounds in lemma~\ref{lem:weightedRV} and its corollaries
(in particular, \eqref{eq:L11powerdecay})
therefore hold on all spaces $L^{p,1}_x(\R^3)$ with $p$ in the appropriate
range.  Similarly, lemma~\ref{lem:SupInverse} asserts that 
$(I+VR_0^\pm(\lambda^2))^{-1}$ is bounded on all $L^{p,1}_x(\R^3)$,
uniformly in $x$.

Two other observations are worth noting at this point.  First is the norm bound
$$\Big\|\frac{V(\cdot)}{|\cdot - x|}\Big\|_{L^{1,1}_x} \les \norm[V][]$$
which holds for all $x \in \R^3$.  
Second, the operator $\frac{d}{d\lambda} \big[e^{-i\lambda|y-x|}
\frac{d}{d\lambda}R_0^+(\lambda^2)\big]$ maps $L^{1,1}_x$ to $L^{\infty,-1}_y$.
This is seen by examining the integration kernel
$$|K(x_2,x_1)| = \big|\frac{d}{d\lambda} e^{i\lambda(|x_2-x_1| - |y-x|)}\big|
\le |x_2 - y| + |x_1 - x|$$
Clearly $\disp \sup_{x_1,x_2} \big|(1+|x_2-y|)^{-1} K(x_2,x_1) (1+|x_1-x|)^{-1}
 \big| \le 2$.

We now return to the remaining term in $\frac{d}{d\lambda}\big[
  e^{-i\lambda|x-y|}B^+(\lambda,x,y)\big]$, namely:
$$\Big\la(I+R_0^+(\lambda^2)V)^{-1}(R_0^+(\lambda^2)V)^k
\Big[\dfrac{d}{d\lambda}\big(e^{-i\lambda|y-x|}\dfrac{d}{d\lambda}
R_0^+(\lambda^2)\big)\Big]
(VR_0^+(\lambda^2))^{m-k}\frac{V(\cdot)e^{i\lambda|\cdot-x|}}{|\cdot-x|},\
\frac{V(\cdot)e^{-i\lambda|\cdot-y|}}{|\cdot-y|}\Big\ra$$
By \eqref{eq:L11powerdecay}, its dual, and the above mapping estimate for the
twice-differentiated resolvent, 
the left-hand function is in $L^{\infty,-1}_y(\R^3)$ 
with norm less than $(1+|\lambda|)^{-(m-6)\eps/4}$.  The right-hand function
is in the dual space $L^{1,1}_y(\R^3)$ with norm controlled by $\norm[V][]$.

The pairing of these two functions is therefore finite, and controlled
pointwise in $\lambda$ by the integrable expression
$(1+|\lambda|)^{-(m-6)\eps/4}$.
\end{proof}

\bibliographystyle{amsplain}

\medskip\noindent
\textsc{Division of Astronomy, Mathematics, and Physics, 253-37 Caltech, Pasadena, CA 91125, U.S.A.}\\
{\em email: }\textsf{\bf mikeg@its.caltech.edu} 

\end{document}